\documentclass[12pt]{amsart}

\usepackage{amsmath, amssymb}
\usepackage{graphicx,color,hyperref}
\usepackage{verbatim, enumitem}

\def\ii{{\rm i}\kern1pt}

\def\proof{\noindent{\it Proof.} }

\def\build#1^#2_#3{\mathrel{\mathop{\null#1}\limits^{#2}_{#3}}}
\def\mertorelbar{\vrule width0.6ex height0.65ex depth-0.55ex}
\def\merto{\mathrel{\mertorelbar\kern1.3pt\mertorelbar\kern1.3pt\mertorelbar
    \kern1.3pt\mertorelbar\kern-1ex\raise0.28ex\hbox{${\scriptscriptstyle>}$}}}

\catcode`\@=11
\newdimen\@rrowlength \@rrowlength=6ex
\def\ssrelbar{\vrule width\@rrowlength height0.64ex depth-0.56ex\kern-4pt}
\def\llra#1{\@rrowlength=#1\ssrelbar\rightarrow}
\catcode`\@=12

\def\semidirect{\mathop{\kern2pt\vrule depth-0.3pt height4.3pt 
\kern-2pt\times}\nolimits}

\newdimen\plainitemindent \plainitemindent=18pt
\def\plainitem#1{\par\noindent
\hangindent\plainitemindent\hbox to\plainitemindent{#1\hss}\ignorespaces}

\catcode`\@=11
\def\openup{\afterassignment\@penup\dimen@=}
\def\@penup{\advance\lineskip\dimen@
  \advance\baselineskip\dimen@
  \advance\lineskiplimit\dimen@}
\newdimen\jot \jot=3pt
\newskip\plaincentering \plaincentering=0pt plus 1000pt minus 1000pt
\def\ialign{\everycr{}\tabskip\z@skip\halign}
\def\eqalign#1{\null\,\vcenter{\openup\jot\m@th
  \ialign{\strut\hfil$\displaystyle{##}$&$\displaystyle{{}##}$\hfil
      \crcr#1\crcr}}\,}
\newif\ifdt@p
\def\displ@y{\global\dt@ptrue\openup\jot\m@th
  \everycr{\noalign{\ifdt@p \global\dt@pfalse \ifdim\prevdepth>-1000\p@
      \vskip-\lineskiplimit \vskip\normallineskiplimit \fi
      \else \penalty\interdisplaylinepenalty \fi}}}
\def\@lign{\tabskip\z@skip\everycr{}} 
\def\displaylines#1{\displ@y \tabskip\z@skip
  \halign{\hbox to\displaywidth{$\@lign\hfil\displaystyle##\hfil$}\crcr
    #1\crcr}}
\def\eqalignno#1{\displ@y \tabskip\plaincentering
  \halign to\displaywidth{\hfil$\@lign\displaystyle{##}$\tabskip\z@skip
    &$\@lign\displaystyle{{}##}$\hfil\tabskip\plaincentering
    &\llap{$\@lign##$}\tabskip\z@skip\crcr
    #1\crcr}}
\def\leqalignno#1{\displ@y \tabskip\plaincentering
  \halign to\displaywidth{\hfil$\@lign\displaystyle{##}$\tabskip\z@skip
    &$\@lign\displaystyle{{}##}$\hfil\tabskip\plaincentering
    &\kern-\displaywidth\rlap{$\@lign##$}\tabskip\displaywidth\crcr
    #1\crcr}}
\def\plaincases#1{\left\{\,\vcenter{\normalbaselines\m@th
    \ialign{$##\hfil$&\quad##\hfil\crcr#1\crcr}}\right.}
\def\plainmatrix#1{\null\,\vcenter{\normalbaselines\m@th
    \ialign{\hfil$##$\hfil&&\quad\hfil$##$\hfil\crcr
      \mathstrut\crcr\noalign{\kern-\baselineskip}
      #1\crcr\mathstrut\crcr\noalign{\kern-\baselineskip}}}\,}

\catcode`\@=12

\def\dlraw{\mathrel{\rlap{$\longrightarrow$}\kern-1pt\longrightarrow}}
\def\vlra{\mathrel{\smash-}\joinrel\mathrel{\smash-}\joinrel%
\kern-2pt\longrightarrow}
\def\srelbar{\vrule width0.6ex height0.65ex depth-0.55ex}
\def\merto{\mathrel{\srelbar\kern1.3pt\srelbar\kern1.3pt\srelbar
    \kern1.3pt\srelbar\kern-1ex\raise0.28ex\hbox{${\scriptscriptstyle>}$}}}

\newdimen\claimskip \claimskip=7pt
\long\def\claim#1|#2\endclaim
{\removelastskip\vskip\claimskip\noindent{\bf#1.}
{\it\ignorespaces#2}\vskip\claimskip\noindent}

\font\ninerm=cmr9
\font\ninei=cmmi9
\font\ninesy=cmsy9
\font\ninebf=cmbx9

\font\nineit=cmti9
\font\ninesl=cmsl9

\font\eightrm=cmr8
\font\eighti=cmmi8
\font\eightsy=cmsy8
\font\eightbf=cmbx8

\font\eightit=cmti8
\font\eightsl=cmsl8
\font\sixrm=cmr6
\font\sixi=cmmi6
\font\sixsy=cmsy6
\font\sixbf=cmbx6
\font\fiverm=cmr5
\font\fivei=cmmi5
\font\fivesy=cmsy5
\font\fivebf=cmbx5

\newfam\itfam
\newfam\slfam
\newfam\bffam
\def\eightpoint{\def\rm{\fam0\eightrm}%
\textfont0=\eightrm \scriptfont0=\sixrm \scriptscriptfont0=\fiverm
 \textfont1=\eighti \scriptfont1=\sixi \scriptscriptfont1=\fivei
 \textfont2=\eightsy \scriptfont2=\sixsy \scriptscriptfont2=\fivesy
 \def\it{\fam\itfam\eightit}%
 \textfont\itfam=\eightit
 \def\sl{\fam\slfam\eightsl}%
 \textfont\slfam=\eightsl
 \def\bf{\fam\bffam\eightbf}%
 \textfont\bffam=\eightbf \scriptfont\bffam=\sixbf
 \scriptscriptfont\bffam=\fivebf
 \normalbaselineskip=9pt
 \setbox\strutbox=\hbox{\vrule height7pt depth2pt width0pt}%
 \normalbaselines\rm}

\def\ninepoint{\def\rm{\fam0\ninerm}%
\textfont0=\ninerm \scriptfont0=\sixrm \scriptscriptfont0=\fiverm
 \textfont1=\ninei \scriptfont1=\sixi \scriptscriptfont1=\fivei
 \textfont2=\ninesy \scriptfont2=\sixsy \scriptscriptfont2=\fivesy
 \def\it{\fam\itfam\nineit}%
 \textfont\itfam=\nineit
 \def\sl{\fam\slfam\ninesl}%
 \textfont\slfam=\ninesl
 \def\bf{\fam\bffam\ninebf}%
 \textfont\bffam=\ninebf \scriptfont\bffam=\sixbf
 \scriptscriptfont\bffam=\fivebf
 \normalbaselineskip=11pt
 \setbox\strutbox=\hbox{\vrule height7pt depth2pt width0pt}%
 \normalbaselines\rm}

\def\plainsection#1{\par\vskip .5cm\penalty -100 
\vbox{\noindent{\sc #1}
\vskip 5pt}
\penalty 500}

\def\Bibitem#1&#2&#3&#4&%
{\hangindent=1.66cm\hangafter=1
\noindent\rlap{\hbox{\rm #1}}\kern1.66cm{\rm #2}{\it #3}{\rm #4.}
\vskip5pt}

\def\square{{\hfill \hbox{
\vrule height 1.453ex  width 0.093ex  depth 0ex
\vrule height 1.5ex  width 1.3ex  depth -1.407ex\kern-0.1ex
\vrule height 1.453ex  width 0.093ex  depth 0ex\kern-1.35ex
\vrule height 0.093ex  width 1.3ex  depth 0ex}}}
\def\bigsquare{{\kern-0.3ex\hbox{
\vrule height 1.7ex  width 0.093ex  depth 0ex\kern-0.093ex
\vrule height 1.8ex  width 1.7ex  depth -1.707ex\kern-0.093ex
\vrule height 1.7ex  width 0.093ex  depth 0ex\kern-1.65ex
\vrule height 0.093ex  width 1.6ex  depth 0ex}\kern0.3ex}}
\def\qed{\phantom{~}$\square$\medskip}
\def\smallskip{\vskip 3pt}
\def\medskip{\vskip 5pt}

\title{On compact pseudoconcave sets}

\author{Zbigniew Slodkowski}
\date{October 2022}

\begin{document}


\begin{abstract} Replying to three questions posed by N. Shcherbina, we show that a compact psudoconcave set can have the core smaller than itself, that the core of a compact set must be pseudoconcave, and that it can be decomposed into compact pseudoconcave sets on which all smooth plurisubharmonic functions are constant.
\end{abstract}

\maketitle 

\plainsection{0. Introduction} 

This paper is inspired by the recent work of N. Shcherbina [Sh55], where he asks under what condition on a compact subset $K$ of a complex manifold $M$ there is a smooth strictly plurisubharmonic function defined in some neighborhood of $K$. In the spirit of [ST04] and [HST17], he defines a subset $c(K)$ of $K$, which (if nonempty) constitutes an obstacle to existence of such a function. This subset, called the {\it core} of K, consists of all points of $p$ of $K$, such that no smooth plurisubharmonic function on $K$ can be strictly plurisubharmonic at $p$. The main result of Shcherbina in [Sh21] is that a strictly plurisubharmonic function on $K$ exists (ie. $c(K)$ is nonempty), if and only if $K$ does not contain any compact pseudoconcave set (cf. Def. 1.1). However, he does not prove this by studying the core,but by examining another subset of $K$, the {\it nucleus}, $n(K)$, which turns out to be the maximal pseudoconcave subset of $K$. 

Shcherbina poses three questions concerning the core. Question 5.1 [Sh21]: "Is $n(K) \subset c(K)$?" We give an example showing it is not true in general (Ex. 2.3). This circumstance makes pertinent his Question 5.2.:"Is the core of a compact set always pseudoconcave?" We show this is true, and we also show, answering his Question 5.3, that the core decomposes into pairwise disjoint pseudoconcave sets on which every smooth plurisubharmonic function is constant (in analogy with similar results in [PS19], and [Sl19]). These facts constitute our main result, Th. 2.1. It might be of some interest that with these properties of the core established, an alternative proof of Shcherbina's main theorem follows easily, cf. Cor. 2.3. in Section 1 we collect some facts about compact pseudoconcave sets that we need.

\claim 0.1. Some conventions in terminology and notation|\rm Following Shcherbina, we say that a funcion is smooth plurisubharmonic on set $K$, if it has a smooth plurisubharmonic extension to some neighborhood of $K$. By {\it smooth} we mean any smoothness class $C^k$, with $k = 0, 1, ...,\infty$. All arguments in Sections 1 and 2 are are correct for all these classes. (In case of $k = 0,1$, {\it strictly} plurisubharmonic should be understood as {\it strongly} plurisubharmonic, in the sens of Richberg [Ri68]). Although nucleus is clearly independent on the smoothness class under consideration, this is unknown for the core (and seems to me unlikely). So the notation $c^k(K)$ might be more accurate, but, like Shcherbina, we suppress the $k$.

A {\it neighborhood} will always be an open set of $M$.

$M$ will denote a complex manifold.

For a subset $A$ of $M$, $bA$ denotes $\overline\setminus A$.      
\endclaim

\claim  Acknowledgement|\rm I am grateful to Nikolai Shcherbina for attracting my attention to these problems, and helpful correspondence concerning Question 5.1.
\endclaim 

\plainsection{1. Background on pseudoconcave sets and Liouville sets}

\claim 1.1. Definition|(a) Let Let $X$ be a locally closed subset of M (ie.
$X\cap(\overline{X}\setminus X) = \emptyset)$. We say that $X$ is a local maximum set, if there do not exist: a point $p\in X$, a neighborhood $V$ of $p$ in $M$, and a $C^{\infty}$-smooth strictly plurisubharmonic function $v$ on $V$, such that $v|V\cap X$ has strict maximum at $p$.

(b) If a local maximum set is closed in $M$, we call it pseudoconcave. 

\endclaim

Below we give just several facts about pseudoconcave sets, that we use in this article. For more background we refer first of all to [Sh21](since  we continue that paper), and then to [Ro55], [Sl86],[ST04],[HST17, 20, 21], [PS19],and [Sl19].

Although pseudoconcavity is the principal concept, the more {\it relaxed} notion of a local maximum set maximum set ads occasionally some flexibility to arguments that would be cumbersome otherwise. See Example 2.3, where the following two properties are used. (They follow immediately from the definition.)

\claim 1.2. Proposition|(a) The union of a family of local maximum subsets of $M$ is pseudoconcave, provided it is closed.

 (b) If $Z \subset X \subset M$, $Z$ relatively closed in $X$, and $X$ a local maximum subset of M, then $X\setminus Z$ is a local maximum subset of M as well.

\endclaim

\claim 1.3. Definition [HST20]| We say that a closed subset $X$ of $M$ is a 
{\it Liouville} set, if every smooth plurisubharmonic function $u$ on $X$, such that $u|X$ is bounded from the above, must be constant on $X$.

\endclaim

\claim 1.4. Proposition|(a) Every Liouville set is pseudoconcave.

 (b) If $K\subset M$ is compact, and $L$ is a compact Liouvile subset of $M$, then 
 $L\subset c(K)$. In particular, $c(L) = L$.

\endclaim

\noindent{\it Comments on proof.} For (a) see [HST20, Lemma 4.1].

 (b) Suppose, to the contrary, that there is a smooth plurisubharmonic function $u$ on a neighborhood $U$ of $X$, which is strictly plurisubharmonic at some neighborhood $V$ of $p\in X$. We can assume $\overline{V}\subset U$ and $X\setminus V\neq\emptyset$. Choose a nonnegative $C^{\infty}$ function  $\rho$ on $U$, such that $\rho(p) > 0$  and $\overline{supp(\rho)}\subset V$. By the standard argument, for some positive $\epsilon$, the function $v:= u + \epsilon\rho$ is plurisubharmonic on $U$. Since $X$ is a compact Liouville set, both $u|X$, and $v|X$ must be constant, which is not possible (as they differ by $\epsilon\rho)$.
 
 \qed
 
 We will say that $X$ is a {\it minimal} pseudoconcave set, if it does not contain any proper pseudoconcave subset $Z$. Bellow we will consider only compact minimal pseudoconcave sets. (A pseudoconcave subset of a compact set is automatically compact.)
 
 \claim 1.5. Proposition| Every compact pseudoconcave subset of manifold $M$ contains a minimal compact pseudoconcave subset.

\endclaim{}

This is very similar to [Sl19, Cor.1.11], and has practically the same proof, but we need this formally stronger statement. (It will be only a consequence of Th.2.1 that the two facts are actually equivalent.)

\proof{} Let $K\subset M$ be a given compact pseudoconcave set. We will define sets $L_t$ by transfinite induction. Let $L_0 = K$. If sets $L_s$ up to ordinal $t$ have already been defined, we let $L_{t+1}$ to be any proper compact pseudoconcave subset of $L_t$, if such exists; if not the induction stops and $L_t$ is the answer. If $t$ is a limit ordinal, we let
$$
L_t:=\bigcap\{L_s: s<t\}.
$$
Clearly $L_t$ is compact and nonempty; we have to check it is a local maximum set. Suppose not, and let $p, V, v: V\leftarrow R$ be as in Def.1.1(i). Choose $B:= B(p,r), r > 0$, a ball in a coordinate neighborhood of $p$, and such that $\overline{B}\subset V$.Let
$F:=\{x\in\overline{B} : v(x)\geq v(p)\}$. Since $F$ is compact, and
$$
\{p\} = F\cap L_t = F\cap\bigcap_{s<t}L_s ,
$$
there is $s_0 < t$, such that $F\cap L_{s_0}\subset B(p,r/2)$.Thus
$$
\max v|L_{s_0}\cap\overline{B(p,r/2)} \geq v(p) > \max v|L_{s_0}\cap bB(p,r).
$$
But this violates another, equivalent, form of local maximum property, namely (ii) in [Sl86, Prop.2.3], so $L_{s_0}$ would not be pseudoconcave. This contradiction shows that $L_t$
is pseudoconcave. For cardinality reasons, the induction has to stop at some ordinal $t$, yielding a compact minimal pseudoconcave subset of $K$.

\qed

\claim 1.6. Proposition| Let $Z$ be a pseudoconcave subset of $M$ and let $u$ be a plurisubharmonic function on $Z$. Suppose $u$ is bounded from the above on $Z$ and attains maximum value at some point of $Z$. Then the (nonempty) set $Y:=\{x\in Z : u(x) = \sup u|Z\}$ is pseudoconcave.

\endclaim{}

\claim|\rm Note $u$ is does not have to be smooth, just uppersemicontinuous. The other assumptions about it hold automatically when $Z$ is compact.

\endclaim

\proof Since $u$ is uppersemicontinuous, $Y$ is closed. Assume WLOG that $\max u|Z = 0$, and  suppose that Y is not a local maximum set. By Def.1.1 there are: a point $p\in Y$, a neighborhood $V$ of $p$, a smooth plurisubfarmonic function $v$ on $V$, and a coordinate ball $B(p,r), r > 0$, such that : $\overline{V}\cap Z$ is compact, $\overline{B(p,r)}\subset V$, $v(p) = 0$ , and 
$v|Y\cap bB(p,r) \leq -1$. So there is a neighborhood $W$ of $Y\cap bB(p,r)$, such that 
$v|W\leq -1/2$. Let $F:= Z\cap bB(p,r)\setminus W $ . Since $F$ is disjoint from $Y$ , $u$ is negative on $F$. Let $c_1:= \max u|F$, then $c_1 < 0$ (as $F$ is compact and $u$ usc). Similarly, $v$ is bounded on $F$ by some finite $c$ .

Choose now $n\in N$ large enough so that $c+nc_1 < - 1/2$, and define a plurisubharmonic function on $V$ by $\phi:= v + nu$.  Now
$\phi\leq v\leq -1/2 $ on $W\cap bB(p,r)$, but 
$\phi\leq c+nc_1 < - 1/2$ on $F$ . 
Thus $\phi|Z\cap bB(p,r) < - 1/2$ , while $\phi(p) = 0$. This contradicts (one of the forms) of local maximum property of $Z$, cf.[Sl86, Prop.2.3(ii)].
 
\qed

\claim 1.7. Corollary|Every minimal compact pseudoconcave set $Z\subset M$ has the Liouville property (with respect to all plurisubharmonic functions on $Z$. 

\endclaim{}

\proof Let $u$ be a plurisubharmonic function defined in a neighborhood of $Z$. As $Z$ is compact, $u|Z$ has maximum value $c$. If $u$ were not constant on Z, then the set 
$Y:=\{x\in Z : u(x)=c \} $ would be different from $Z$. Since, by last proposition, $Y$ is pseudoconcave, we would obtain contradiction with the minimality of $Z$.

\qed

We Recall the properties of the {\it regularized maximum function}.

\claim 1.8. Proposition.[HM88,4.13-14]|For every $\delta > 0$ there is a smooth convex function $M_{\delta}$ on $R^2$, with both partial nonnegative partial derivatives, and such that 
$$
M_{\delta}(t, s) = \max(t,s) \hbox{ , when }  |t-s|\geq \delta .
$$

\endclaim

\plainsection{2. Pseudoconcave decomposition of the core}

\claim 2.1. Theorem| Let $K\subset M$ be compact, and $p\in K$.  Let $Y$ be the set of all points $x\in K$ such that $u(x) = u(p) $ for all smooth plurisubharmonic functions  $u$ on $K$.

\emph{(a)}  If $Y = \{p\} $ , then then $p$ does not belong to the core $c(K)$.

\emph{(b)} If $Y$ is not a singleton, then $Y$ is a (compact) Liouville set. In particular, $Y$ is pseudoconcave, and $c(Y) = Y$.

\endclaim

\proof{}Suppose that we are given a neighborhood $V$ of $Y$, and a smooth plurisubharmonic function $v:V\rightarrow [1/2, 1]$. During most of the proof we will construct, for the pair $(V, v)$, a smooth plurisubharmonic function $\chi^*$ on $K$, with special properties (cf. Conclusion below). Then a pair $(V, v)$ will be chosen differently in cases (a) and (b), and using corresponding  $\chi^*$, we will obtain (a) and (b) immediately.

By the definition of $Y$, there is a family of smooth plurisubharmonic functions 
$\phi_t: U_t\rightarrow R$, with $U_t$ neighborhoods of $K$, for $t\in T$, such that 
$$
Y = \bigcap_{t\in T}\{x\in K : \phi_t (x) = 0\}.
$$
Since the sets in this intersection are compact, we can choose finitely many of them so that (simplifying the notation slightly) we get
$$
Y\subset Z:=\bigcap_{i=1}^n\{x\in K : \phi_i (x) = 0\}\subset V,
$$
with $\phi_i:U_i\rightarrow R$,  $K\subset U_i$, $i = 1,...,n$. Choose a relatively compact open set $U$, and such that
$$
K\subset U\subset\overline{U}\subset U_0:=\bigcap_{i=1}^n U_i.
$$
Shrinking $U_i$'s if needed, we we can assume WLOG that all $\phi_i$'s are uniformly bounded on $U_0$.

In this part of the proof we will follow the beginning of the proof of Lemma 3.5 in[Sl19] (preserving the notation). Define first two smooth plurisubharmonic functions $\phi$ and $\mu$ on $U_0$. Let
$$
\phi(x):= \phi_1(x) +...+\phi_n(x) \hbox{ , for }  x\in U_0.
$$
To define $\mu$, let first 
$v_{\epsilon}(t_1,...,t_n):= t_1+...+t_n+\epsilon(t_1^2+...+t_n^2)$ which defines a smooth convex function on $R^n$, for $\epsilon> 0$. Since $m_0:= \inf\{\phi_i(x): x\in U_0, i=1,...,n\}$ is finite (by the choice of $U_0$), we can choose an $\epsilon> 0$, such that
$1+\epsilon m_0 > 0$ . Then $\frac{\partial v_{\epsilon}}{\partial t_i}(t) > 0$ when 
$\min(t_1,...,t_n) > m_0$, in particular on the joint range of the vector valued function
$(\phi_1,...,\phi_n)$. So, by Lemma 1.13(iii) in [Sl19], the function 
$\mu(x):= v_{\epsilon}(\phi_1(x),...,\phi_n(x))$ is smooth plurisubharmonic on $U_0$. Similarly as in [Sl19], we have $\phi(x)\leq\mu(x)$, for $x\in U_0$, and
$$
\{x\in U_0 : \phi(x) = \mu(x)\} = \bigcap_{i=1}^n\{x\in U_0 : \phi_i(x) = 0\}.
$$

Observe that
$$
Z = K\cap\bigcap_{c>0}\{x\in\overline{U} : |\phi_i(x)|\leq c, i=1,...,n\}.
$$
Since the  sets $ \{x\in\overline{U} : |\phi_i(x)|\leq c, i=1,...,n\}$ are compact and decreasing with $c$, there is $c > 0$, such that $K\cap\overline{H}\subset V$, where
$H:=\bigcap_{i=1}^n\{x\in\overline{U} : |\phi_i(x)|< c\}$. Finally, there is an open set $W$ such that
$$
K\subset W \subset\overline{W}\subset U \hbox{ , and   }          \overline{W}\cap\overline{H}\subset V. 
$$ 
Observe that
$$
U\setminus H\subset \bigcup_{i+1}^n\{x\in U : |\phi_i(x)| = 1 \}\subset \{x\in U : \mu(x) > \phi(x)\}.
$$
Now, as $W\setminus H$ is relatively compact in $U$, we conclude that
$\Delta:= \inf\{ (\mu-\phi)(x): x\in W\setminus H \} > 0$. So, choosing a $\delta\in (0,\Delta/5)$, we obtain
$$
\mu\geq\phi +\Delta\geq(\phi+(\Delta/2)v) +\delta  \hbox{ , on } W\cap V\setminus H.
$$
By this inequality and Prop.1.8, the function $\chi$ defined as $M_{\delta}(\phi+(\Delta/2)v, \mu)$ on 
$W\cap V$, equals $\mu$ on $W\cap V\setminus\overline{H}$. Function $\chi$ is smmoth plurisubharmonic on $W\cap V$. Define now function $\chi^*$ as $\chi$ on $V\cap W$ and as $\mu$ on $W\setminus\overline{H}$.
Since the latter two functions are equal on the intersection of their domains 
$(V\cap W)\cap(W\setminus\overline{H}) = W\cap V\setminus\overline{H}$, function $\chi^*$ is well defined smooth plurisubharmonic function on the union of these domains:
$$
(V\cap W)\cup(W\setminus\overline{H})= W \hbox{ , since } W\cap\overline{H}\subset V. 
$$

Consider now $\chi^*$ near $Y$. Since $\phi|Y = \mu|Y =0$, it follows that $\phi+(\Delta/2)v > \mu + \delta$, on $Y$, and so on some neighborhood of $Y$, thus $\chi^* = \phi+(\Delta/2)v$ near $Y$.
 
\claim Conclusion| For a smooth plurisubharmonic function $1/2<v< 1$ on a neighborhood of $Y$, there is a smooth plurisubharmonic function $\chi^*$ on $K$ that is equal to $\phi+(\Delta/2)v$ near $Y$, and to $(\Delta/2)v$ on $Y$.

\endclaim

\noindent {\it Proof of {\rm (a).}} If $Y = \{p\}$, we choose as $V$ a small coordinate ball centered at $P$, and as $v$ a smooth strictly plurisubharmonic function on $V$, bounded by $1/2$ and $1$. By the Conclusion, we obtain a smooth plurisubharmonic function $\chi^*$ on $K$, that is strictly plurisubharmonic near $p$. Thus
$p$ is not in $c(K)$.

\noindent {\it Proof of {\rm (b).}}Suppose $Y$ is not a Liouville set. Then there is a smooth plurisubharmonic function $v$ on some open $V$ containing $Y$, and $p,q\in Y$, such that $v(p)\neq v(q)$ (WLOG $1/2<v<1$). By the Conclusion, there is a smooth strictly plurisubharmonic function $\chi^*$ on $K$, that is equal to $(\Delta/2)v$ on $Y$. But then $\chi^*(p)\neq\chi^*(q)$, which contradicts the definition of $Y$.

\qed

Observe that the last theorem gives positive answer to Shcherbina's Question5.3. Also, since $c(K)$ is the union of sets $Y$ having local maximum property by (b), $c(K)$ is pseudoconcave by Proposition 1.2(a). Finally, Theorem 2.1 allows for an alternative proof of Shcherbina's Main Theorem.

\claim 2.2. Corollary| Let $K$ be a compact subset of a complex manifold. Then there is a smooth strictly plurisubharmonic function on $K$ , if and only if $K$ does not contain any compact pseudoconcave subset.

\endclaim

\proof{ If $K$ contains a compact pseudoconcave set $X$, then, by by Proposition 1.5, $X$ contains a minimal compact pseudoconcave set $Z$, which by Proposition 1.7 is a Liouville set, and so, by Proposition 1.4(b), must be containd in the core $c(K)$. Core being nonempty, a smooth strictly plurisubharmonic fuction on $K$ cannot exist.

If $K$ does not contain any compact pseudoconcave set, then $c(K)$ must be empty (it would be otherwise pseudoconcave), which means that a smooth strictly plurisubharmonic fuction on $K$ does exist.
}
\qed

The next example illustrates Theorem 2.1, and gives answer to Shcherbina's Question 5.1.

\claim 2.3. Example|\rm Let $M:= \overline{C}\times C $ and 
$K:=\overline{C}\times S^1 \cup\{0\}\times\overline{D(0, 2)} $. Consider pluriharmonic functions
$u(z, w)=\Re w$, and $v(z, w)=\Im w$,  and let $Y_{a,b}:=\{(z, w)\in K : u(z, w)=a,  v(z, w)= b \}$. 
If $a^2 + b^2 =1$, then $Y_{a,b}=\overline{C}\times \{a+ib\}$ and is a compact Liouville set. Otherwise
$Y_{a,b} = \{(0, a+ib)\}$ is a single point (when $a^2 + b^2 < 1$, or $1 < a^2 + b^2\leq 4$). By Theorem 2.1 the core is the union of Liouville sets contained in $K$. Hence $c(K) =\overline{C}\times S^1 $. 

Let now $K_1:= c(K)\cup\{0\}\times D(0, 1)$. Since this set is closed, and is the union of local maximum sets (the second not closed), it is pseudoconcave by Proposition 1.2(a), and so is contained in $n(K)$. Thus in this case the nucleus and the core are not equal (which answers [Sh21,Question 5.1]). Observe yet that 
$n(K) = K_1$. Suppose $n(K)$ contains $K_1$ properly. By Proposition 1.2(b), $n(K)\setminus K_1$ would be a relatively closed local maximum subset of the annulus $\{0\}\times \{a+ib : 1< a^2+b^2\leq 4 \}$. But this is not possible. 
            
\endclaim

\plainsection{References}

\parskip 1.5pt plus 0.5pt minus 0.5pt

\Bibitem[HST17]&Harz, T., Shcherbina, N., Tomassini, G.:& On defining functions and cores for unbounded
domains I.& Math. Z. {\bf 286} (2017), 987-1002&

\Bibitem[HST20]&Harz, T., Shcherbina, N., Tomassini, G.:& On defining functions and cores for unbounded
domains II.& J. Geom. Anal. {\bf 30} (2020), 2293-2325&

\Bibitem[HST21]&Harz, T., Shcherbina, N., Tomassini, G.:& On defining functions and cores for unbounded
domains III.& Mat. Sb. {\bf 212}, 6 (2021), 126-156&

\Bibitem[HL88]&Henkin,G.-M ., Leiterer, J., G.:& Andreootti-Grauert theory by integral formulas.& Progress in Mathematics {\bf 74} (1988)&

\Bibitem [MST18]&Mongodi, S., Slodkowski, Z., Tomassini, G.:& Weakly complete complex surfaces.& Indiana  Univ. Math. J. {\bf 67} (2018), 899-935&

\Bibitem[PS19]&Poletsky, E. A., Shcherbina, N.:& Plurisubharmonically separable complex manifolds. & Proc. Amer. Math. Soc. {\bf 147} (2019), 2413 - 2424&

\Bibitem[Ri68]&Richberg, R.:& Stetige streng pseudokonvexe Funktionen.& Math. Ann. {\bf 175} (1968), 257--286&

\Bibitem[Ro55]&Rothstein, W.:& Zur Theorie der analytischen Manningfaltigkeiten in Raume von n komplexen Veranderlichen.& Math. Ann. {\bf 129} (1955), 96-138&

\Bibitem[Sh21]&Shcherbina, N.V.:& On compacts possessing strictly plurisubharmonic functions. & Izvestiya : Mathematics  {\bf 85} (2021), 605-618&

\Bibitem[ST04]&Slodkowski, Z., Tomassini, G.:& Minimal kernels of weakly complete spaces. & J. Funct. Anal. {\bf 210} (2004), 125--147&

\Bibitem[Sl19]&Slodkowski, Z.:&Pseudoconcave decompositions in complex manifolds.& Contemp. Math. {\bf 735} (2019), 239-259&

\vskip3mm\noindent
Zbigniew Slodkowski
Department of Mathematics, University of Illinois at Chicago\\
851 South Morgan Street, Chicago, IL 60607, USA\\
\emph{e-mail}\/: zbigniew@uic.edu

\end{document}